\documentclass[12pt]{article}

\newtheorem{cor}{Corollary}[section]
\newtheorem{prop}{Proposition}[section]
\newtheorem{lem}{Lemma}[section]
\newtheorem{Def}{Definition}[section]

\newcommand{\be}{\begin{equation}}

\newcommand{\ee}{\end{equation}}
\font\BBb=msbm10 at 12pt
\newcommand{\Bbb}[1]{\mbox{\BBb #1}}
\newcommand{\qed}{\hspace*{\fill}Q.E.D.}
\title{\bf Manifolds with Quadratic Curvature Decay and Slow Volume Growth}
\author{John Lott \thanks{Supported by NSF Grant DMS-9704633.} \\
Department of Mathematics\\
University of Michigan\\
Ann Arbor, MI  48109-1109
\and
Zhongmin Shen \\
Department of Mathematics\\
Indiana University-Purdue University at Indianapolis\\
Indianapolis, IN 46202-3216}  
\date{September 17, 1998}
\begin{document}
\maketitle

\begin{center}
{\it To Detlef Gromoll on his 60th birthday}
\end{center}

\begin{abstract}
We show that there are topological obstructions for a noncompact 
manifold to
admit a Riemannian metric with quadratic curvature decay and a volume growth
which is slower than that of the Euclidean space of the same dimension.
\end{abstract}
\newpage
\section{Introduction}

A major theme is Riemannian geometry is the relationship between curvature
and topology.  For compact manifolds, one can constrain the curvature and
diameter and ask whether one obtains topological restrictions on the
manifold.  If the manifold is noncompact then a replacement for a diameter
bound is a constraint on how the curvature behaves in terms of the distance
from a basepoint.  More precisely,
let $M$ be a complete connected $n$-dimensional Riemannian manifold.
Fix a basepoint $m_0 \in M$.

\begin{Def}
\label{defn2} $M$ has quadratic curvature decay (with constant $C>0$)
if  for all $m\in M$ and all $2$-planes $P$ in $T_mM$, the sectional
curvature $K(P)$ of $P$ satisfies 
\begin{equation}
|K(P)|\le C/d(m_0,m)^2.  \label{eqn1}
\end{equation}
\end{Def}

Note that condition (\ref{eqn1}) is scale-invariant in that
it is unchanged by a constant rescaling of the Riemannian metric.
Many interesting results have been obtained when
the sectional curvature has some faster-than-quadratic curvature decay
\cite{Abresch (1985),Greene-Petersen-Zhu (1994)}.
In this paper we concentrate instead on the case of quadratic curvature decay.
In itself, condition (\ref{eqn1}) does not put any restrictions on the
topology of a manifold.  One can show that any connected smooth paracompact
manifold has a Riemannian metric with quadratic curvature decay;
see \cite[p. 96]{Gromov} or Lemma \ref{lemmaquad} below. On the other hand,
we will show that if in addition one restricts the volume growth of the
metric, then one does obtain topological restrictions on $M$. The first
question is whether $M$ has finite topological type.

\begin{Def}
\label{defn1} $M$ has finite topological type if 
$M$ is homotopy-equivalent to a
finite $CW$-complex.
\end{Def}

\begin{Def}
$M$ has lower quadratic curvature decay (with constant $C > 0$)
if  for all $m\in M$ and all $2$-planes $P$ in $T_mM$, the sectional
curvature $K(P)$ of $P$ satisfies 
\begin{equation}
K(P) \ge -  C/d(m_0,m)^2.  \label{eqn1a}
\end{equation}
\end{Def}

Let $B_t$ denote the metric ball of radius $t$ around $m_0$ and let $S_t$
denote the distance sphere of radius $t$ around $m_0$. 
If $M$ has lower quadratic curvature decay then
by a standard argument,
one can show that $M$ has at most polynomial volume growth; see 
Lemma \ref{lem2} below.

\begin{prop}\label{prop_order2}
Suppose that $M$ has lower quadratic curvature decay.
If ${\rm vol}(B_t) = o(t^2)$ as $t \rightarrow \infty$
and $M$ does not collapse at infinity, i.e.
$\inf_{x\in M} {\rm vol}(B_1(x)) >0$, then $M$ has finite topological type.
\end{prop}

The $o(t^2)$ bound in Proposition \ref{prop_order2} cannot be improved to
$O(t^2)$, as shown in Example 3 below.

Next, we consider manifolds with volume growth slower than that of the
Euclidean space of the same dimension.

\begin{Def}
\label{defn3} $M$ has slow volume growth if
\begin{equation}
\liminf_{t\rightarrow\infty }{\rm vol}(B_t)/t^n\ =0.  \label{eqn2}
\end{equation}
\end{Def}

There is a notion of an {\em end} $E$ of $M$ and of $E$ being
contained in an open set ${\cal O} \subset M$; see, for example, 
\cite[p. 80]{Bonahon}.

\begin{Def}
\label{defn4}
An end $E$ of $M$ is {\em tame} if it is contained in an open set
diffeomorphic to $(0, \infty) \times X$ for some smooth 
connected closed manifold $X$.
\end{Def}

We remark that $X$ is determined by $E$ only up to $h$-cobordism.
Hereafter we assume that $M$ is oriented.

\begin{prop}
\label{prop1}
Suppose that $M$ has quadratic curvature decay and slow volume growth.
Let $E$ be a tame end of $M$ as in Definition \ref{defn4}.
Then for any
product $\prod_k p_{i_k}(TX)$ of
Pontryagin classes of $X$ and
any bounded cohomology class $\omega \in {\rm H}^{l}
(X; {\Bbb R})$ with $l + 4\sum_k i_k = n-1$,
\begin{equation}
\int_X \omega \cup \prod_k p_{i_k}(TX) = 0.
\end{equation}
\end{prop}

\noindent
\textbf{Example :} There is no metric of quadratic curvature decay and slow
volume growth on $\mbox{\BBb R } \times \Bbb C {\rm P}^{2k}$.\\

Next, we give a sufficient condition for $M$ to have a
metric of quadratic
curvature decay and slow volume growth.

\begin{prop}
\label{prop2} Let $X$ be a closed manifold with a polarized $F$-structure 
\cite{Cheeger-Gromov (1986)}. Suppose that $X = \partial N$ for some 
smooth compact
manifold $N$. Then there is a complete Riemannian metric on $M = {\rm Int}(N)$
of quadratic curvature decay and slow volume growth.
\end{prop}

It follows from Proposition \ref{prop2} that when $n$ is even, there is a
metric on ${\Bbb R}^n$
of quadratic curvature decay and slow volume growth. The case when $n$ is
odd is less obvious.
 
\begin{prop}
\label{prop3} For all $n > 1$, there is a complete
Riemannian metric on ${\Bbb R}^n$ 
of quadratic curvature decay and slow volume growth.
\end{prop}

If $X$ is a closed oriented manifold with a polarized $F$-structure then the
Pontryagin numbers and Euler characteristic of $X$ vanish. Based on
Proposition \ref{prop2}, one may think that under the hypotheses of
Proposition \ref{prop1}, one could also show that the Euler characteristic
of $X$ vanishes. However, Proposition \ref{prop3} shows that this is not the
case, as the Euler characteristic of $S^{n-1}$ is two if $n$ is odd.

We can combine Propositions \ref{prop1}-\ref{prop3} 
to obtain some low-dimensional results.

\begin{cor}
\label{cor1}
Let $N$ be a smooth compact connected oriented manifold-with-boundary
of dimension $n$.\\
1. If $n = 2$ then ${\rm Int}(N)$ has a metric of quadratic curvature decay
and slow volume growth.\\
2. If $n = 3$ then ${\rm Int}(N)$ has a metric of quadratic curvature decay
and slow volume growth if and only if $\partial N$ consists of $2$-spheres
and $2$-tori.\\
2. If $n = 4$, suppose that Thurston's Geometrization Conjecture holds.
Then ${\rm Int}(N)$ has a metric of quadratic curvature decay
and slow volume growth if and only if the connected components of $\partial N$
are graph manifolds.
\end{cor}

Finally, as in \cite{Cheeger-Gromov (1985)}, there is an integrality result
for the integral of the Gauss-Bonnet-Chern form,
which we state without proof.

\begin{prop}
\label{prop4}
Suppose that $M$ has a complete Riemannian metric $g$
of quadratic curvature decay with ${\rm vol}(B_t) = o(t^n)$ and
$\int_1^\infty \frac{{\rm vol}(B_t)}{t^n} \: \frac{dt}{t} \: < \: \infty$.
Let $e(M, g) \in \Omega^n(M)$ be the Gauss-Bonnet-Chern form.
Then $\int_M e(M, g) \in {\Bbb Z}$.
\end{prop}

We thank M. Gromov for pointing out the relevance of bounded cohomology.

\section{Examples}
\noindent
1. Let $N$ be a smooth compact connected 
$n$-dimensional manifold-with-boundary. Let $h$ be a
metric on $\partial N$. Given $c \ge 1$, consider the metric on 
$[1, \infty) \times \partial N$ given by 
$dt^2 \: + \: t^{2c} \: h$. Extend this to a smooth metric $g$
on ${\rm Int}(N) = N \cup_{\partial N} [1, \infty) \times \partial N$.
Then $g$ has quadratic curvature decay and polynomial volume growth.  
By choosing $c$ large,
the degree of  volume growth can be made arbitrarily large. Taking $c = 1$,
we see that having quadratic curvature decay and volume growth of order
$O(t^n)$ in no way restricts the topology of the ends.
\\ \\
2. For $c \in {\Bbb R}$, consider the metric on $[1, \infty) \times S^1$ given
by $dt^2 \: + \: t^{2c} \: d\theta^2$. 
Cap this off by a disk at $\{1\} \times 
S^{1}$ to obtain a smooth metric $g$ on ${\Bbb R}^2$.  Then $g$ has quadratic
curvature decay. If $c \: < \: -1$ then $({\Bbb R}^2, g)$ has finite volume.
Hence the assumption of quadratic curvature decay gives no nontrivial
lower bound on volume growth.\\ \\
3. Start with the Euclidean metric on the annulus 
$A = \{ (x,y) \in {\Bbb R}^2 : 1 \le x^2 + y^2 \le 4 \}$.
Add a handle to ${\rm Int}(A)$, keeping the metric the same near
$\partial A$. Let $h$ denote the corresponding metric
on $T^2 - D^2$. With an obvious notation, for $j \in {\Bbb N}$, let
$2^j \cdot (T^2 - D^2)$ denote the rescaled metric. 
Consider 
$(T^2 - D^2) \cup_{S^1} 2 \cdot (T^2 - D^2) \cup_{S^1} 4 \cdot (T^2 - D^2)
\cup_{S^1} \ldots$ with its corresponding metric.
Cap it off with a disk to obtain a smooth metric $g_\Sigma$ 
on an infinite genus surface
$\Sigma$. For $n \ge 2$,
let $g_{T^{n-2}}$ be a flat metric on the $(n-2)$-torus. Then
the product metric $\left( \Sigma, g_\Sigma \right) \times \left(T^{n-2},
g_{T^{n-2}} \right)$ has quadratic curvature decay, volume growth of
order $t^2$ and infinite topological type.  This shows that the $o(t^2)$
condition in Proposition \ref{prop_order2} cannot be improved to $O(t^2)$.

\begin{lem} \label{lemmaquad}
If $M$ is a smooth connected paracompact manifold then $M$ admits a complete
Riemannian metric of quadratic curvature decay.
\end{lem}
{\it Proof}:
First, $M$ admits a complete Riemannian metric $h$
of bounded sectional curvature \cite{Greene (1978)}.
Given $\phi \in C^\infty \left( M \right)$, put $g =
e^{2\phi} h$. We have
\begin{equation}
R^i_{\: jkl}(g) = R^i_{\: jkl}(h) \: - \: \widetilde{\phi}^i_{\: k} \:
h_{jl} \: + \: \widetilde{\phi}^i_{\: l} \: h_{jk} \: - \: \delta^i_{\: k}
\: \widetilde{\phi}_{jl} \: + \: \delta^i_{\: l} \: \widetilde{\phi}_{jk} \:
- \: \phi_{;r} \: \phi^{;r} \: \left( \delta^i_{\: k} \: h_{jl} \: - \:
\delta^i_{\: l} \: h_{jk} \right),  \label{eqn17}
\end{equation}
where $\widetilde{\phi}_{ab} \: = \: \phi_{;ab} \: - \: \phi_{;a} \:
\phi_{;b}$. Let $d_h$ denote the distance function with respect to $h$ and
let $d_g$ denote the distance function with respect to $g$.
By \cite[Theorem 1.8]{Cheeger-Gromov (1991)}, there is a 
$\phi \in C^\infty(M)$ and a constant $c > 0$ such that\\
1. $\phi(m) \: \le \: d_h(m_0, m) \: \le \: \phi(m) \: + \: c$.\\
2. $\parallel \nabla \phi \parallel_\infty \: \le \: c$.\\
3. $\parallel {\rm Hess}(\phi) \parallel_\infty \: \le \: c$.\\
Then from (\ref{eqn17}), in order to show that $g$ has quadratic curvature 
decay it suffices to show that there is a constant $C > 0$ such that
$d_g(m_0, m) \: \le \: C \: e^{\phi(m)}$ for all $m \in M$. 
Let $\gamma$ be a normalized minimal
geodesic, with respect to $h$, from $m_0$ to $m$.  Then measuring the
length of $\gamma$ with respect to $g$,
\begin{equation}
d_g(m_0,m) \le \int_0^{d_h(m_0,m)} e^{\phi(\gamma(t))} dt \le
\int_0^{d_h(m_0,m)} e^{t} dt = e^{d_h(m_0,m)} - 1 \le
e^c \: e^{\phi(m)}.  
\end{equation}
The lemma follows. \qed

\section{Proof of Proposition \ref{prop_order2}}
First of all, we show that every manifold with lower
quadratic Ricci curvature decay has polynomial volume growth.

\begin{lem} \label{lem2}
Suppose that there is a constant $C > 0$ such that
for each $m\in M$ and each unit
vector $v\in T_mM$, the Ricci curvature satisfies
\begin{equation}
{\rm Ric} (v, v) \geq -(n-1) {C\over d(m_0, m)^2}, \label{Ricdecay}
\end{equation}
Put $N=(n-1)\: \frac{\sqrt{1+4C}-1}{2} \:  + \: n$.
Then there is a constant  $C_0=C_0(n, C)>0$
such that for $t\geq 3$,
\begin{equation} {\rm vol} (B_t) \leq   C_0 {\rm vol}(S_1) \: t^{N} + 
{\rm vol}(B_1)\label{comparison1}
\end{equation}
and
\begin{equation}
 {\rm vol} (B_{t+1}- B_{t-1})
\leq C_0\frac{{\rm vol}(B_{t-1})}{t-1}.\label{comparison2}
\end{equation}
\end{lem}
{\it Proof}:
Let $\Pi_t = {1\over n-1} \sum_{i=1}^{n-1}k_i$ denote the mean curvature 
of the regular part of $S_t$,
where $\{k_i\}_{i=1}^{n-1}$ are the principal curvatures.
Letting
$dA_t$ and $dA_{m_0}$
denote the volume forms on $S_t$ and $S_{m_0}M$ respectively,
define $\varphi_t: S_{m_0}M \to S_t$ by
\begin{equation} \varphi_t (v) = \exp_{m_0}(tv)
\end{equation}
and  define
$\eta_t: S_{m_0}M \to (0, \infty)$ by
\begin{equation}
 (\varphi_t)^* dA_t |_v= \eta_t(v) \: dA_{m_0}.
\end{equation}
We have
\begin{equation}
{\rm vol}(S_t) =\int_{S_{m_0}M} \eta_t (v) \: dA_{m_0}
\end{equation}
and
\begin{equation}
(n-1) \: \Pi_t|_{\varphi_t(v)} = \eta_t'(v)/\eta_t(v). \label{Eta}
\end{equation}
As $t \rightarrow 0$,
\begin{equation}
(n-1) \: \Pi_t|_{\varphi_t(v)} = {{n-1}\over t} - {{\rm Ric} ( v,v) \over 3} 
t + o(t).\label{Taylor}\end{equation}
Put $\Pi(t) = \Pi_t|_{\varphi_t(v)}$
and $v(t) = (\exp_{m_0})_*(tv)$. The Riccati equation 
implies
\begin{equation}
\Pi'(t) + \Pi(t)^2 \leq -{ {\rm Ric}(v(t), v(t)) \over n-1} . \label{Riccati}
\end{equation}
Put $\alpha = \frac{ \sqrt{1+4C} +1}{2}$ and consider
\begin{equation}
 f(t) = e^{\int_1^t \Pi(s)(s)ds} 
\Big [ t^{\alpha} \Pi(t) - \alpha t^{\alpha -1} \Big ].
\end{equation}
Then (\ref{Taylor}) implies that $\lim_{t\to 0+} f(t) =0$.
On the other hand, 
from (\ref{Ricdecay}) and (\ref{Riccati}), we have
\begin{equation}
f'(t) = t^\alpha \: e^{\int_1^t \Pi(s)(s)ds}\Big [\Pi'(t)
+ \Pi(t)^2 - {\alpha (\alpha -1) }t^{-2} \Big ]
\leq 0.
\end{equation}
Thus $f(t) \leq 0$, giving
\begin{equation}
\Pi(t) \leq \alpha t^{-1} .
\end{equation}
 Together with (\ref{Eta}), we conclude that
$ \eta_t(v)/t^{(n-1)\alpha}$ is nonincreasing. 
  This implies
that ${\rm vol}(S_t)/t^{(n-1)\alpha}$ is nonincreasing, too.
 As
\begin{equation}
{\rm vol}(B_t)-{\rm vol}(B_1) 
= \int_1^t \frac{{\rm vol}(S_s)}{s^{(n-1)\alpha}} \: s^{(n-1)\alpha} \: ds,
\end{equation}
we obtain
\begin{equation}
{\rm vol}(S_1) \: \int_1^t s^{(n-1)\alpha} \: ds \: \ge \: {\rm vol}(B_t)- {\rm vol}(B_1)
 \: \ge \:
\frac{{\rm vol}(S_t)}{t^{(n-1)\alpha}} \: \int_1^t s^{(n-1)\alpha} \: ds.
\end{equation}
Hence
\begin{equation} {\rm vol} (B_t) \:
\leq \: \frac{1}{(n-1)\alpha + 1} 
\: {\rm vol}(S_1)\: t^{(n-1)\alpha + 1 } + {\rm vol}(B_1).
\end{equation}
Also,
\begin{eqnarray}
 {\rm vol}(B_{t+1}-B_{t-1}) & = & \int_{t-1}^{t+1}
\frac{{\rm vol}(S_s)}{s^{(n-1)\alpha}} \: s^{(n-1)\alpha} \: ds \\
&\le &
\frac{{\rm vol}(S_{t-1})}{(t-1)^{(n-1)\alpha} }\int_{t-1}^{t+1}  \: s^{(n-1)\alpha} \: ds
\nonumber \\
&\leq &  \frac{{\rm vol}(B_{t-1})-{\rm vol}(B_1)}{ \int_{1}^{t-1} s^{(n-1)\alpha}ds }
\int_{t-1}^{t+1} s^{(n-1)\alpha} \: ds \nonumber \\
& \leq & C_0 \: \frac{ {\rm vol}(B_{t-1}) }{ t-1} \nonumber
\end{eqnarray}
for large enough $C_0$. \qed \\ \\
{\bf Proof of Proposition \ref{prop_order2}}:\\
We use critical point theory of the distance function; for a review, see
\cite{Cheeger}.
Let us say that a connected component $\Sigma_t$ of $S_t$ is {\em good} if 
it is part of the  boundary of an unbounded component of $M-{B}_t$ and
there is a ray from $m_0$ passing through $\Sigma_t$.
\begin{lem}
Suppose that there is a $t_0 > 0$ such that 
if $t \ge t_0$ then there is no critical point of $d_{m_0}$ on any good
component $\Sigma_t$ of $S_t$. Then $M$ has finite topological type.
\end{lem}
{\it Proof}:
Let $E$ be an end of $M$. We know that there is a normalized ray
$\gamma$ such that  $\gamma(0) = m_0$ and $\gamma$
exits $E$. By assumption, for all $t \ge t_0$, the connected component
$\Sigma_t$ of
$S_t$ which contains $\gamma(t)$ does not include any critical points
of $d_{m_0}$. By the isotopy lemma \cite[Lemma 1.4 and p. 35]{Cheeger}, 
it follows
that the unbounded component $U$ of $M - B_{t_0}$ containing
$\{\gamma(t)\}_{t = t_0}^\infty$ is homeomorphic to $[0, \infty) \times
\Sigma_{t_0}$, with $\Sigma_{t_0}$ a closed connected topological manifold. 
In particular, $U \cap S_t$ is connected and good, so $U$ does not
contain any critical points.
As $S_{t_0}$ is compact, it has a finite number of
connected components. It follows that $M-B_{t_0}$ has a finite number of
bounded connected components and a finite number of unbounded 
connected components. Thus there is some $t_1 > t_0$ such that
$M - B_{t_1}$ does not have any critical points, from which the lemma follows.
In fact, the proof shows that
$M$ is homeomorphic to the interior of a compact topological
manifold-with-boundary. \qed \\

Define
\begin{equation} {\cal D}(m_0, t) = \sup {\rm Diam}(\Sigma_t),
\end{equation}
where the supremum is taken over all good components $\Sigma_t$ of $S_t$
and the diameter is measured using the metric on $M$. 
We claim that 
if the manifold has lower quadratic curvature decay and if
\begin{equation}
 \lim_{t\to \infty} { {\cal D}(m_0, t) \over t} =0 
\end{equation}
there is a $t_0 > 0$ such that 
if $t \ge t_0$ then there is no critical point of $d_{m_0}$ on any good
component $\Sigma_t$ of $S_t$.

For a pair of points $p, q\in M$, define 
$$ e_{pq}(x) = d(p, x)+d(q, x)-d(p, q).$$
Clearly, for any $t > 0$ and
any point $m\in M-B_{2t}$ on a ray from $m_0$ which intersects
$\Sigma_t$,  
\begin{equation}
e_{m_0 m} (x) \leq 2 {\cal D}(m_0, t) \: \: \: {\rm for} \: x\in \Sigma_t.
\end{equation}
By assumption, the sectional curvature on $M-B_{t/2}$ satisfies
\begin{equation} K_M \: \geq \: - \: \frac{4C}{t^2}.
\end{equation}
Assume that there is a $t_0 > 0$ such that for $t >t_0$,
\begin{equation}
 {\cal D}(m_0, t) \: \leq \: \frac{t}{4\lambda \sqrt{C}}, \label{diam}
\end{equation}
where 
$\lambda$ is a large constant which will be specified later.

Suppose that $x\in \Sigma_t$ is a critical point of $d_{m_0}$.
Take a minimizing geodesic $\tau$ from $x$ to $
m$. There is a minimizing geodesic $\sigma$ from
$x$ to $m_0$ such that
$\angle (\dot{\sigma}(0), \dot{\tau}(0)) \le {\pi \over 2}$.
Take two points $p=\sigma(a)$ and $q=\tau(a)$ where
$a = \frac{t}{\lambda\sqrt{C}}$. By the triangle inequality, we have
\begin{equation} 
e_{pq}(x)\leq e_{m_0m}(x) \leq 2 {\cal D}(m_0, t).
\end{equation}
For $\lambda\geq \frac{100}{\sqrt{C}}$, we see that the triangle $\Delta_{pxq}$
is contained in a small neighborhood of $x$ inside $M-B_{t/2}$. Then 
we can apply the Toponogov inequality to $\Delta_{pxq}$ and 
obtain 
\begin{equation} \cosh (c_0 d(p, q)) \:
\leq \: \cosh^2 (c_0 a),
\end{equation}
where $c_0 = \frac{2\sqrt{C}}{t}$.
Note that
\begin{equation}
 c_0 d(p, q) = c_0[ 2a- e_{pq}(x)]\ge 2c_0 [a- {\cal D}(m_0, t)] \geq 
\frac{3}{\lambda}.
\end{equation} 
We obtain
\begin{equation} \cosh \left( \frac{3}{\lambda} \right) \leq 
\cosh^2 \left( \frac{2}{\lambda} \right).
\end{equation}
This is impossible for sufficiently large $\lambda$. 

\bigskip
Finally, we must show that if ${\rm vol}(B_t) = o(t^2)$
and if there is a $v > 0$ such that ${\rm vol}(B_1(x)) > v$ for all $x \in M$,
then (\ref{diam}) holds for large $t$. 

Let $ \Sigma_t $ be a connected component of the boundary of
an unbounded component of $M- B_t$.
For any $x, y\in \Sigma_t$, there is
a continuous curve $c:[0, r]\to \Sigma_t$ from $x$ to $y$. Suppose that
$d(x,y) > 2$.
Then there is a partition
$0=t_0 < t_1 < \cdots < t_k= r$ such that
$\{B_{1}( c(t_i))\}_{i=0}^k$ are disjoint and 
$ B_2 (c(t_i)) \cap B_2(c(t_{i+1}))\not=\emptyset$.
Note that $B_1(c(t_i)) \subset B_{t+1} -\overline{B_{t-1}}$.
We have
\begin{equation} (k+1) v \leq \sum_{i=0}^k  {\rm vol} (B_1(c(t_i)))
\leq {\rm vol} ( B_{t+1} -\overline{B_{t-1}} )
\leq C_0 \frac{ {\rm vol} ( B_{t-1})}{t-1}.\end{equation}
Thus
\begin{equation} 
{\rm Diam} (\Sigma_t) \leq \sum_{i=0}^{k-1} d(c(t_i), c(t_{i+1})) 
\leq C_1 \frac{ {\rm vol} ( B_{t-1})}{t-1},\end{equation}
giving
\begin{equation} 
\lim_{t\to \infty} {{\cal D}(m_0, t) \over t} =0 .
\end{equation}
This proves Proposition \ref{prop_order2}. \qed

\section{Proof of Proposition \ref{prop1}}

Fix an open set ${\cal O}$ containing $E$
which is diffeomorphic to $(0, \infty) \times X$. For 
$u > 1$, let $\widehat{M}$ denote $M$ with the metric $u^{-2} \: g_M$. Let 
$\widehat{{\cal O}}$ denote the copy of ${\cal O}$ in $\widehat{M}$. 
Let $\widehat{B}_t$ and 
$\widehat{S}_t$ denote the metric ball and metric sphere in $\widehat{M}$
around $m_0$. Rescaling (\ref{eqn1}), there is a constant $C^\prime > 0$
such that the region $\widehat{B}_{100} - \widehat{B}_{1/100}$ has sectional
curvatures bounded by $C^\prime$, uniformly in $u$. Put 
\begin{equation}
T_{1/10} (\widehat{S}_1 \cap \widehat{{\cal O}}) = 
\{ \widehat{m} \in \widehat{M}
\: : \: d(\widehat{m}, \widehat{S}_1 \cap \widehat{{\cal O}}) \le 1/10 \}.
\label{eqn3}
\end{equation}
By \cite[Theorem 0.1]{Cheeger-Gromov (1991)}, there is a constant $C^{\prime
\prime} > 0$ independent of $u$ such that there is a connected codimension-$0$
submanifold $U_u$ of $\widehat{M}$ with 
\begin{equation}
(\widehat{S}_1 \cap \widehat{{\cal O}}) 
\subset U_u \subset T_{1/10} (\widehat{S}_1
\cap \widehat{{\cal O}}),  \label{eqn4}
\end{equation}
\begin{equation}
{\rm vol}(\partial U_u) \le C^{\prime \prime} \: 
{\rm vol}( T_{1/10} (\widehat{S}_1
\cap \widehat{{\cal O}}))  \label{eqn5}
\end{equation}
and 
\begin{equation}
\parallel \Pi_{\partial U_u} \parallel \: \le \: C^{\prime \prime},
\label{eqn6}
\end{equation}
where $\Pi_{\partial U_u}$ is the second fundamental form of $\partial U_u$
in $\widehat{M}$. Then by the Gauss-Codazzi equation, the intrinsic
sectional curvature of $\partial U_u$ is uniformly bounded in $u$. Rescaling
to $M$, we have 
\begin{equation}
{\rm vol}( T_{1/10} (\widehat{S}_1 \cap \widehat{E})) = 
u^{-n} \: {\rm vol}( T_{u/10}
(S_u \cap {\cal O})) \le u^{-n} \: {\rm vol}( B_{11u/10}).  \label{eqn7}
\end{equation}

Let $\{u_j\}_{j=1}^\infty$ be a sequence in ${\Bbb R}^+$ approaching
infinity such that 
\begin{equation}
\lim_{j \rightarrow \infty} {\rm vol} \left( B_{11u_j/10} \right)/u_j^n = 0.  \label{eqn8}
\end{equation}
For $j$ large, let $Y_j$ be a connected component of $\partial U_{u_j}$.
Let ${\cal O}_j$ be the oriented cobordism between $Y_j$ and $X$ coming from
the unbounded component of $M - Y_j$ corresponding to $E$, truncated at
some level $\{R_j\} \times X$.  Let
$i : Y_j \rightarrow {\cal O}_j$ be the inclusion and let $\pi : {\cal O}_j
\rightarrow (0, \infty) \times X \rightarrow X$ be projection.
Then
\begin{eqnarray}
&& \int_X \omega \cup \prod_k p_{i_k}(TX) - 
\int_{Y_j} (\pi \circ i)^* \omega \cup \prod_k p_{i_k}(TY_j) = \\
&&\int_{{\cal O}_j} d \left( \pi^* \omega \wedge \prod_k
p_{i_k}(T{\cal O}_j) \right) = 0. \nonumber
\end{eqnarray}
From (\ref{eqn5}), (\ref{eqn7}), (\ref{eqn8}) and 
\cite[p. 37]{Gromov},
we have that 
$\int_{Y_j} (\pi \circ i)^* \omega \cup \prod_k p_{i_k}(TY_j) = 0$
if $j$ is large enough. The proposition follows. \qed

\section{Proof of Proposition \ref{prop2}}

Suppose that $\{g(t)\}_{t \in [1, \infty)}$ is a smooth $1$-parameter family
of Riemannian metrics on $X$ with sectional curvatures that are uniformly
bounded in $t$. Then one can check that $dt^2 + t^2 \: g(t)$ is a metric of
quadratic curvature decay on $[1, \infty) \times X$ if $\parallel g^{-1}(t)
\: \frac{dg}{dt} \parallel_\infty = O \left( \frac{1}{t} \right)$
and $\parallel g^{-1}(t) \: \frac{d^2 g}{dt^2} \parallel_\infty = O \left( 
\frac{1}{t^2} \right)$. Put $\delta = t^{-1}$ and let $g(t)$ be the
Riemannian metric on $X$ defined in \cite[Section 3]{Cheeger-Gromov (1986)}.
Then $\{g(t)\}_{t \in [1, \infty)}$ has uniformly bounded sectional curvature
in $t$. We claim that $\parallel g^{-1}(t) \: \frac{dg}{dt} \parallel_\infty
= O \left( \frac{1}{t} \right)$ and $\parallel g^{-1}(t) \: 
\frac{d^2 g}{dt^2} \parallel_\infty = O \left( \frac{1}{t^2} \right)$.
The metric $g(t)$ is defined by a finite recursive process. One starts with
an invariant Riemannian metric $g_0$ for the $F$-structure and puts $g_1(t) =
\log^2 (1 + t) \: g_0$. Clearly $\parallel g_1^{-1}(t) \: \frac{dg_1}{dt}
\parallel_\infty = O \left( \frac{1}{t} \right)$ and $\parallel
g_1^{-1}(t) \: \frac{d^2 g_1}{dt^2} \parallel_\infty = O \left( \frac{1}{t^2}
\right)$. Then 
\begin{equation}
g_{j+1}(t) = 
\left\{
\begin{array}{ll}
\rho_j^2 \: g_j^{\prime}(t) \: + \:
h_j(t) & \mbox{ on $U_j$}, \\
g_j(t) & \mbox{ on $X - U_j$.}
\end{array}
\right.
\label{eqn9}
\end{equation}
where\\
1. $U_j$ is a certain open subset of $X$,\\
2. $g_j^\prime(t)$ is the
part of $g_j(t)$ 
corresponding to tangent vectors to the $F$-structure on $U_j$, \\
3. $h_j(t)$ is the part of $g_j(t)$ corresponding to normal vectors to
the $F$-structure on $U_j$ and \\
4. $\rho_j = t^{- \: \frac{\log(f_j)}{\log(1/2)}}$
with $f_j : X \rightarrow [1/2,1]$ a
certain smooth function which is identically one on $X - U_j$.

It follows by induction on $j$ that there is a metric of quadratic curvature
decay and small volume growth on $[1, \infty) \times X$. Gluing $[1, \infty)
\times X$ onto $N$, we obtain the desired metric on $M$. \qed

\section{Proof of Proposition \ref{prop3}}

If $n$ is even then $S^{n-1}$ has a polarized $F$-structure coming from a free 
$S^1$-action and the result follows from Proposition \ref{prop2}. The first
nontrivial case is when $n = 3$.

Suppose that $n = 3$. By \cite[Example 1.4]{Cheeger-Gromov (1985)}, there is
a metric $h$ on ${\Bbb R}^3$ with finite volume and bounded sectional
curvature. Our metric will be conformally related to $h$. Let us first give
the construction of $h$ in detail. For $j \in {\Bbb Z}^+$, let $C_j$ be the
complement of a small solid torus in a solid torus. Then topologically, 
\begin{equation}
{\Bbb R}^3 = (S^1 \times D^2) \cup_{T^2} C_1 \cup_{T^2} C_2 \cup_{T^2}
\ldots  \label{eqn10}
\end{equation}
We take $m_0 \in S^1 \times D^2$. Each $C_j$ can be decomposed as $C_j =
(\Sigma_{2j} \times S^1_{2j}) \cup_{T^2} (\Sigma_{2j+1} \times S^1_{2j+1})$,
where $\Sigma_{2j}$ is a $2$-sphere with three disks removed, $\Sigma_{2j+1}$
is a $2$-disk and $S^1_{2j}$, $S^1_{2j+1}$ are circles. In \cite[Figure 1.3]
{Cheeger-Gromov (1985)}, $\Sigma_{2j}$ is represented as a rectangle with a
disk removed and with the vertical sides identified. Put $\partial
\Sigma_{2j} = S^1_{2j,1} \cup S^1_{2j,2} \cup S^1_{2j,3}$, where $S^1_{2j,1}$
is the top side of the rectangle, $S^1_{2j,2}$ is the bottom side of the
rectangle and $S^1_{2j,3}$ is the circle enclosing the removed disk. Put 
$\partial \Sigma_{2j+1} = S^1_{2j+1, 1}$. The identifications of the toroidal
boundaries are 
\begin{eqnarray}  \label{eqn11}
S^1_{2j+1,1} \times S^1_{2j+1} & \sim &S^1_{2j,2} \times S^1_{2j}, \\
S^1_{2j,3} \times S^1_{2j} & \sim & S^1_{2j-2,1} \times S^1_{2j-2}, 
\nonumber
\end{eqnarray}
where 
\begin{eqnarray}  \label{eqn12}
S^1_{2j+1,1} & \sim & S^1_{2j}, \\
S^1_{2j+1} & \sim & S^1_{2j,2},  \nonumber \\
S^1_{2j,3} & \sim & S^1_{2j-2},  \nonumber \\
S^1_{2j} & \sim & S^1_{2j-2,1}.  \nonumber
\end{eqnarray}

We will put product metrics on $\Sigma_{2j} \times S^1_{2j}$ and 
$\Sigma_{2j+1} \times S^1_{2j+1}$. Let $\epsilon_i$ be the length of $S^1_i$
and let $\delta_{i,*}$ be the length of $S^1_{i,*}$. Then (\ref{eqn12}) gives
the relations 
\begin{eqnarray}  \label{eqn13}
\delta_{2j,1} & =& \epsilon_{2j+2}, \\
\delta_{2j,2} & =& \epsilon_{2j+1},  \nonumber \\
\delta_{2j,3} & = &\epsilon_{2j-2},  \nonumber \\
\delta_{2j+1,1} & =& \epsilon_{2j}.  \nonumber
\end{eqnarray}

We will take $\epsilon_i = e^{-i}$. Let $\Sigma_\infty$ be a thrice-punctured
sphere with a Riemannian metric such that three ends $E_1, E_2, E_3 \cong
(1, \infty) \times S^1$ are isometric to $dr^2 \: + \: e^{-2r} \: d\theta^2$.
Put $\Sigma_0 = \Sigma_\infty - (E_1 \cup E_2 \cup E_3)$. Let $u \in
C^\infty([0,1])$ be a nondecreasing function such that 
\begin{equation}
\left\{
\begin{array}{ll}
u(s) = s & \mbox{ if $s \in [0,\frac{1}{3}]$}, \\
1 & \mbox{ if $s \in [\frac{1}{2},1]$.}
\end{array}
\right.
\end{equation}
Given $k \in {\Bbb Z}^+$, put $E(k) = [0, k] \times S^1$
with the metric $dr^2 \: +
\: e^{- 2k u(r/k)} \: d\theta^2$. Then put 
\begin{equation}
\Sigma_{2j} = \Sigma_0 \cup_{\partial \Sigma_0} (E(2j+2) \cup E(2j+1) \cup
E(2j-2)),  \label{eqn15}
\end{equation}
isometrically. Similarly, let $\Sigma^\prime_\infty$ be a once-punctured
sphere with a Riemannian metric such that the end $E \cong (1, \infty)
\times S^1$ is isometric to $dr^2 \: + \: e^{-2r} \: d\theta^2$. Put 
$\Sigma_0^\prime = \Sigma_\infty^\prime - E$ and 
\begin{equation}
\Sigma_{2j+1} = \Sigma_0^\prime \cup_{S^1} E(2j),  \label{eqn16}
\end{equation}
isometrically. Then one can check that $\{\Sigma_i\}_{i=1}^\infty$ have
uniformly bounded volume and curvature. Glue together the product metrics on 
$\{\Sigma_{2j} \times S^1_{2j}\}_{j=1}^\infty$ and $\{\Sigma_{2j+1} \times
S^1_{2j+1}\}_{j=1}^\infty$ to give the metric $h$ on ${\Bbb R}^3$. As 
$\sum_{j=1}^\infty e^{-j} < \infty$, it follows that $h$ has bounded curvature
and finite volume.

Given $\phi \in C^\infty \left({\Bbb R}^3 \right)$, put $g =
e^{2\phi} h$. 
By (\ref{eqn17}), the weighted sectional curvatures 
\begin{equation}
\left\{e^{2\phi(m)} \: |K(P,g)| \right\}_{m \in M, \: P \subset
T_mM}  \label{eqn18}
\end{equation}
are uniformly bounded provided that the gradient $\nabla \phi$ of $\phi$ and
the Hessian $H(\phi)$ of $\phi$ are uniformly bounded with respect to $h$.

We construct $\phi$ on $\Sigma_{2j} \times S^1_{2j}$ and $\Sigma_{2j+1}
\times S^1_{2j+1}$ to be the pullbacks of functions on $\Sigma_{2j}$ and 
$\Sigma_{2j+1}$, respectively. Let $\phi_\infty \in C^\infty(\Sigma_\infty)$
be a Morse function with one critical point, of saddle type, such that 
\begin{eqnarray}  \label{eqn19}
\phi_\infty \big|_{E_1} & = &40 \: d(\cdot, \Sigma_0), \\
\phi_\infty \big|_{E_2} & = &10 \: d(\cdot, \Sigma_0),  \nonumber \\
\phi_\infty \big|_{E_3} & = &- \: 80 \: - \: 40 \: d(\cdot, \Sigma_0), 
\nonumber \\
\phi_\infty(\Sigma_0) & \subset & [-80, 0].  \nonumber
\end{eqnarray}
Then in terms of (\ref{eqn15}), put 
\begin{equation}
\phi \big|_{\Sigma_{2j}} = 80j^2 + 80j + \phi_\infty \big|_{\Sigma_{2j}}.
\label{eqn20}
\end{equation}
Similarly, let $\phi^\prime_\infty \in C^\infty(\Sigma^\prime_\infty)$ be a
Morse function with one critical point, a local maximum, such that 
\begin{eqnarray}  \label{eqn21}
\phi^\prime_\infty \big|_{E} & = & - \: 10 \: d(\cdot, \Sigma^\prime_0), \\
\phi^\prime_\infty (\Sigma^\prime_0) & \subset & [0,10].  \nonumber
\end{eqnarray}
Then in terms of (\ref{eqn16}), put 
\begin{equation}
\phi \big|_{\Sigma_{2j+1}} = 80j^2 + 120j + 10 + \phi^\prime_\infty 
\big|_{\Sigma_{2j+1}}.  \label{eqn22}
\end{equation}
Finally, define $\phi$ on the $S^1 \times D^2$ factor in (\ref{eqn10}) so as
to extend $\phi$ to a smooth function on ${\Bbb R}^3$.

It is easy to see that $\nabla \phi$ and $H(\phi)$ are uniformly bounded on 
${\Bbb R}^3$. As 
\begin{equation}
d_g(m_0,m)^2 \: |K(P,g)| = \frac{d_g(m_0,m)^2}{e^{2\phi(m)}} \: e^{2\phi(m)}
\: |K(P,g)|,  \label{eqn23}
\end{equation}
in order to show that $g$ has quadratic curvature decay, it suffices to show
that $e^{- \phi(m)} \: d_g(m_0,m)$ is uniformly bounded with respect to $m
\in {\Bbb R}^3$. Let $T^2$ be the first torus factor in (\ref{eqn10}).
Then it suffices to show that $e^{- \phi(m)} \: d_g(T^2,m)$ is uniformly
bounded with respect to $m \in {\Bbb R}^3$. Let $\{\gamma(s)\}_{s \in
[0,t]}$ be a piecewise smooth path from $m$ to $T^2$ which is unit-speed
with respect to $h$, and along which $\phi$ is nonincreasing. Then letting 
$L_g(\gamma)$ denote the length of $\gamma$ with respect to $g$, we have 
\begin{equation}
e^{- \phi(m)} \: d_g(T^2,m) \: \le \: e^{- \phi(m)} \: L_g(\gamma) =
\int_0^t e^{\phi(\gamma(s)) - \phi(m)} \: ds.  \label{eqn24}
\end{equation}
We take $\gamma$ to be (reparametrized) gradient flow of $\phi$ starting
from $m$. Although $\phi$ is not a Morse function, we note that gradient
flow on $\Sigma_{2j} \times S^1_{2j}$ is essentially the same as gradient
flow on $\Sigma_{2j}$, as it is constant in the $S^1_{2j}$-factor, and
gradient flow on $\Sigma_{2j+1} \times S^1_{2j+1}$ is essentially the same
as gradient flow on $\Sigma_{2j+1}$, as it is constant in the 
$S^1_{2j+1}$-factor. 
If the projection of $\gamma$ onto $\Sigma_{2j}$ or $\Sigma_{2j+1}$
meets a critical point $c$ of saddlepoint type, we extend $\gamma$ beyond $c$
to become a piecewise smooth curve with a corner, again following a downward
gradient trajectory. We continue this process until $\gamma$ hits $T^2$.
Changing variable to $u = \phi(m) - \phi(\gamma(s))$, we have 
\begin{equation}
\int_0^t e^{\phi(\gamma(s)) - \phi(m)} \: ds = \int_0^{\phi(m)} e^{-u} \: 
\frac{du}{|\nabla \phi|(\phi^{-1}(u))}.  \label{eqn25}
\end{equation}

As $\phi(\gamma(s))$ is nonincreasing, if $m \in C_j$ then $\gamma$ never
enters $\Sigma_{2k+1} \times S^1_{2k+1}$ for $k < j$. Also $\gamma$ hits at
most one critical point in each $\Sigma_{2k}$ for $k < j$. By the
construction of $\phi$, if $c_k \in \Sigma_{2k}$ is the critical point then 
$\phi \big|_{c_k \times S^1_{2k}} \in [80k^2 + 80 k - 80, 80 k^2 + 80k]$.
Thus the singularities of $\frac{1}{|\nabla \phi|(\phi^{-1}(u))}$ are
well-spaced in $u$. If $\gamma$ passes through a critical point $c$ and $u_0
= \phi(c)$ then $\frac{1}{|\nabla \phi|(\phi^{-1}(u))} \sim 
\frac{1}{\sqrt{|u-u_0|}}$ for $u \sim u_0$. 
From the uniform nature of $\nabla \phi$ near
the critical points, it follows that there is a constant $D > 0$,
independent of $m \in {\Bbb R}^3$, such that for all $x \in [0, \phi(m)
- 1]$, 
\begin{equation}
\int_x^{x+1} \frac{du}{|\nabla \phi|(\phi^{-1}(u))} \le D.  \label{eqn26}
\end{equation}
Then 
\begin{equation}  \label{eqn27}
\int_0^{\phi(m)} e^{-u} \: \frac{du}{|\nabla \phi|(\phi^{-1}(u))} \le 
\frac{D}{1-e^{-1}}.
\end{equation}
Thus $g$ has quadratic curvature decay.

Put $t_{j+1} = d(m_0, C_{j+1})$. For $j > 0$, each path from $m_0$ to 
$C_{j+1}$ must pass through $C_j$. Put 
\begin{equation}  \label{eqn28}
D_j = (S^1 \times D^2) \cup_{T^2} C_1 \cup_{T^2} \ldots \cup_{T^2} C_j.
\end{equation}
Then $B_{t_{j+1}}(m_0) \subset D_j$ and so 
${\rm vol}(B_{t_{j+1}}) \le {\rm vol}(D_j)$. 
With respect to (\ref{eqn15}), let $F_j$ be the subset 
$\left[j+2, 2j+2 \right] \times S^1_{2j} \subset E(2j+2) \times S^1_{2j}$. For
large $j$, $\phi \big|_{D_j - F_j} \le 80j^2 + 120j + 80$ and so 
\begin{equation}  \label{eqn29}
{\rm vol} (D_j - F_j) \: \le \: e^{240j^2 + 360j +240} \: {\rm vol}(\Bbb R^3, h).
\end{equation}
On the other hand, 
\begin{equation}  \label{eqn30}
{\rm vol}(F_j) = \int_{j+2}^{2j+2} e^{3(80j^2 + 80j + 40x)} \: e^{-2(2j+2)} \: dx
\: = \: \frac{1-e^{-120j}}{120} \: e^{240j^2+480j+240} \: e^{-2(2j+2)}.
\end{equation}
Thus 
\begin{equation}  \label{eqn31}
{\rm vol}(B_{t_{j+1}}) = O \left( e^{240j^2+480j+240} \: e^{-2(2j+2)} \right).
\end{equation}
As any path from $m_0$ to $C_{j+1}$ must pass through $F_j$, 
\begin{equation}  \label{eqn32}
t_{j+1} \: \ge \: \int_{j+2}^{2j+2} e^{80j^2 + 80j + 40x} \: dx = \frac{1 -
e^{-40j}}{40} \: e^{80j^2 + 160j + 80}.
\end{equation}
Thus 
\begin{equation}  \label{eqn33}
{\rm vol}(B_{t_{j+1}})/t_{j+1}^3 \: = \: O \left( e^{-2(2j+2)} \right),
\end{equation}
showing that $g$ has slow volume growth.

If $n > 3$, we can do a similar construction in which $C_j$ is the
complement of a small $T^{n-2} \times D^2$ in $T^{n-2} \times D^2$ and $C_j$
is decomposed as $(\Sigma_{2j} \times T^{n-2}) \cup_{T^{n-1}} (\Sigma_{2j+1}
\times T^{n-2})$. \qed

\section{Proof of Corollary \ref{cor1}}
\noindent
1. If $n = 2$, put a metric on ${\rm Int}(N)$ with flat cylindrical
ends.\\
2. If $n = 3$, suppose that $\partial N$ consists of $2$-spheres and
$2$-tori. For a $2$-sphere component of $\partial N$, put a metric
coming from Proposition \ref{prop3} on
the corresponding end of ${\rm Int}(N)$. For a
$2$-torus component of $\partial N$, put a flat metric on the corresponding
end $(1, \infty) \times T^2$ of ${\rm Int}(N)$. 
This gives the desired metric on ${\rm Int}(N)$.
Now suppose that ${\rm Int}(N)$ has a metric with quadratic curvature decay
and slow volume growth.  From Proposition \ref{prop1}, the simplicial volume of
$\partial N$ must vanish. Thus $\partial N$ consists of $2$-spheres and
$2$-tori.\\
3. If $n = 4$, suppose that the connected components of $\partial N$ are graph
manifolds. Then $\partial N$ has a polarized $F$-structure and Proposition
\ref{prop2} implies that there is a metric on ${\rm Int}(N)$ with quadratic
curvature decay and slow volume growth. Now suppose that Thurston's
Geometrization Conjecture holds and that ${\rm Int}(N)$ has a metric with
quadratic curvature decay and slow volume growth. 
From Proposition \ref{prop1}, the simplicial volume of
$\partial N$ must vanish. From \cite{Soma}, this implies that the
connected components of $\partial N$ are graph manifolds. \qed

\end{document}